\documentstyle[12pt,subeqnarray,epsf]{article}

\begin{document}
\newcommand{\dis}{\displaystyle}
\newcommand{\dml}{\rm dim ~ }
\newcommand{\krl}{\rm ker ~ ~ }
\newcommand{\expon}{\rm e}
\newcommand{\id}{ 1 \hspace{-2.85pt} {\rm I} \hspace{2.5mm}}
\newcommand{\func}{ (\frac{\Phi(\Gamma) - \Phi(J_0)}{\Gamma -
[J_0][J_0 + 1]})^{\frac{1}{2}} }
\newcommand{\D}{\Delta}
\newcommand{\nJ}{\tilde{J}}
\newcommand{\q}{{q^\prime}}
\newcommand{\ca}{{\cal C}}
\newcommand{\U}{{\cal U}}
\newcommand{\alp}{\alpha^\prime}
\newcommand{\beq}{\begin{equation}}
\newcommand{\eeq}{\end{equation}}
\newcommand{\bsq}{\begin{subeqnarray}}
\newcommand{\esq}{\end{subeqnarray}}
{\thispagestyle{empty}
\rightline{} 
\rightline{} 
\rightline{} 
\vskip 1cm
\centerline{\large \bf A General }
\centerline{\large \bf $q$-Oscillator Algebra}


\vskip 2cm
\centerline{L. C. Kwek 
{\footnote{E-mail address:
scip6051@leonis.nus.edu.sg }},
and C. H. Oh 
{\footnote{E-mail address:
phyohch@leonis.nus.edu.sg }}
}
\centerline{{\it Department of Physics, Faculty of Science, } }
\centerline{{\it National University of Singapore, Lower Kent Ridge,} }
\centerline{{\it Singapore 119260, Republic of Singapore. } }}
\vskip 0.1in

\vskip 1cm
\centerline{\bf Abstract} \vspace{5mm}

\noindent{
It is well-known that the Macfarlane-Biedenharn $q$-oscillator and its
generalization has no Hopf structure, whereas the Hong Yan $q$-oscillator
can be endowed with a Hopf structure. In this letter, 
we demonstrate that it is
possible to construct a general $q$-oscillator algebra which includes
the Macfarlane-Biedenharn oscillator algebra and the Hong Yan oscillator
algebra as special cases.
}
\newpage


The harmonic oscillator has often emerged as a basic theoretical tool
for investigating a variety of physical systems.  It appears naturally
in quantum optics and the theory of coherent radiation.  The $q$-analog
of the harmonic oscillator was first proposed by Arik and Coon \cite{arik}
and subsequently applied to the realization of the $q$-deformed $su(2)$
algebra by Macfarlane and Biedenharn \cite{mac,bied}. A good
introduction to the $q$-deformed oscillators and its representation
can be found in ref \cite{kulish}

Macfarlane and Biedenharn defined the $q$-deformed harmonic oscillator
algebra or the $q$-boson algebra as that generated by the operators $\{
a, a^\dagger, N \}$ obeying the commutation relations
\begin{subeqnarray}
\lbrack N, a^\dagger \rbrack &=& a^\dagger, \\
\lbrack N, a \rbrack &=& - a, \\
a a^{\dagger} - q a^{\dagger} a & = & q^{-N}.
\label{macf}
\end{subeqnarray}
The casimir operator for the algebra, $\ca_1$, is given by
\beq
\ca_1  =  q^{-N} (a^\dagger a - [N]). \label{casmac1}
\eeq
In the undeformed situation, namely in the limit $q \rightarrow 1$, this
central element becomes trivial since its eigenvalues can be shifted to
zero and by von Neumann theorem, all representations of the oscillators
are unitarily equivalent to each other.  This is not the case for the
$q$-deformed oscillator.  Different values of the central element can 
label different representations.   

The representation theory for a class of $q$-deformed oscillator algebras
defined in terms of arbitrary function of the number operator $N$ has
also been developed by Quesne and Vansteenkiste \cite{quesne}.  More
recently, Irac-Astaud and Rideau \cite{irac1,irac2}
have constructed Bargmann
representations corresponding to these generalized $q$-deformed
oscillator algebras and showed that Bargmann representations can exist
for some deformed harmonic oscillators which admit non-Fock represntations.

The Macfarlane-Biedenharn $q$-deformed oscillator does not seem to
possess a Hopf structure.  Hong Yan \cite{yan} has proposed another
$q$-deformed oscillator that can admit a Hopf structure.  The
commutation relation for this $q$-oscillator are given by
\begin{subeqnarray}
\lbrack N, a^\dagger \rbrack &=& a^\dagger, \\
\lbrack N, a \rbrack &=& - a, \\
\lbrack a, a^{\dagger} \rbrack & = & [N + 1] - [N]  \label{hyo}
\end{subeqnarray}
where the bracket $[x]$ denotes $\dis \frac{q^x - q^{-x}}{q - q^{-1}}$.
The casimir $\ca_2$ for this algebra is given by
\beq
\ca_{2}  =  a^\dagger a - [N].\label{casmac}
\eeq
This $q$-deformed oscillator algebra can be contracted from the
$q$-deformed $su(2)$ algebra by a generalized Inoue-Wigner
transformation \cite{kwek1}.  
Nevertheless if we impose positive norm requirement for the states, then
at the representation level, the identification of Hong Yan
$q$-oscillator algebra with the $q$-deformed $su(2)$ algebra 
can break down for some
values of $|q| =1$.   In fact, the positive norm requirement
\cite{fuji1} is in
conflict with the truncation condition \cite{oh2}
imposed on the states of the
oscillator so as to get finite multiplets for $su_{\sqrt{q}}(2)$.  In
short, for $|q| =1$ ($q= \expon^{i \epsilon}, ~ ~ \epsilon$
arbitrary) Hong Yan oscillator algebra is different from
$su_{\sqrt{q}}(2)$ algebra.
Note that the
Macfarlane-Biedenharn and Hong Yan algebra are distinct, but in the
usual $q$-Fock representation, they become equivalent.  

The Macfarlane-Biedenharn oscillator can be written in several
equivalent forms. For Arik and Coon \cite{arik} $q$-deformed
oscillator, the commutation
relation for the operators $a$ and $a^\dagger$ in eq(\ref{macf}c) is
expressed as
\beq
a a^\dagger - q a^\dagger a = 1
\eeq
For Chaturvedi and Srinivasan \cite{chatur} $q$-oscillator, 
they define the commutation
relation above as
\beq
a a^\dagger - a^\dagger a = q^{-N}
\eeq
and for Chakrabarti and Jagannathan \cite{chakra} $q$-oscillator, 
they modify the commutation relation so that
it accommodates two parameters, $q_1$ and $q_2$, giving the relation
\beq
a a^\dagger - q_1 a^\dagger a = q_2^{-N}.
\eeq

Many other generalizations have been constructed
\cite{quesne, irac1,irac2}. The purpose of this letter is to devise a new
$q$-oscillator algebra such that it not
only reduces to the Hong Yan algebra for a specific choice of parameters
but also becomes the Macfarlane-Biedenharn $q$-oscillator when one of the
two parameters vanishes.  Thus our general $q$-oscillator algebra
straddles across the Macfarlane-Biedenharn algebra and the Hong Yan
algebra. This feature is, to our knowledge, not present in any
previously known generalized $q$-oscillator algebras.  In the following
we first write down the generalized Macfarlane-Biedenharn 
algebra and show that it
can be reduced to previously known forms of the Macfarlane-Biedenharn
algebra. We then proceed to present our general oscillator.

Following Duc \cite{duc}, we may generalize
the Macfarlane-Biedenharn $q$-oscillator
algebra \cite{kwek1} by introducing two additional parameters
$\alpha$ and $\beta$. The generalized commutation relations for the
Macfarlane-Biedenharn algebra are 
\begin{subeqnarray}
\lbrack N, a^\dagger \rbrack &=& a^\dagger, \\
\lbrack N, a \rbrack &=& - a, \\
a a^{\dagger} - q^{\alpha} a^{\dagger} a & = & q^{\beta N}, \\
\ca_3 & = & q^{-\alpha N} (a^\dagger a - [N]_{\alpha,\beta}), \label{genmacf}
\end{subeqnarray}
where $\dis [x]_{\alpha,\beta} = \frac{q^{\alpha x}- q^{\beta x}}
{q^\alpha - q^\beta}$ is a generalized $q$-bracket.  Despite its
complexity, this algebra is not a new one.  It can be reduced 
to the
usual Macfarlane-Biedenharn $q$-oscillator algebra.  To see this, 
one can define new operators $\dis A= q^{-\frac{(\alpha + \beta) N}{4}} a,
A^\dagger = a^\dagger q^{-\frac{(\alpha + \beta) N}{4}} $ and  map
the generalized commutation relations eq(\ref{genmacf}) to
\begin{subeqnarray}
\lbrack N, A^\dagger \rbrack &=& A^\dagger, \\
\lbrack N, A \rbrack &=& - A, \\
A A^{\dagger} -  q^\prime 
A^{\dagger} A & = & q^{\prime -N}. \label{genmacf1}
\end{subeqnarray}
where $\dis q^\prime = q^{\frac{\alpha - \beta}{2}}$.
Clearly, one gets the
form of eq(\ref{macf}) by identifying $q^{\beta -
\alpha}$ as $(q^\prime)^{-2}$.  Indeed, with this generalized form of
the Macfarlane-Biedenharn $q$-oscillator algebra, one can see easily
that some of the well-known $q$-oscillator algebras are essentially the
same as the Macfarlane-Biedenharn algebra and can be mapped from one to
another by invertible transformation.  To be specific, we note that the
Arik-Coon oscillator \cite{arik},  
Chaturvedi-Srinivasan $q$-oscillator \cite{chatur} and the
Chakrabarti-Jaganathan \cite{chakra} oscillator correspond to the cases 
when $\alpha =1, \beta = 0$, $\alpha = 0, \beta = 1$ and $\alpha,\beta$
arbitrary respectively.

We now propose a new general $q$-oscillator algebra 
in which the operators $\{ a, a^\dagger, N \}$ satisfy the following  
relations
\begin{subeqnarray}
\lbrack N, a^\dagger \rbrack &=& a^\dagger, \\
\lbrack N, a \rbrack &=& - a, \\
\lbrack a, a^{\dagger} \rbrack & = & [N + 1]_{\alpha,\beta} - 
[N]_{\alpha,\beta} \\
\ca_{4} & = & q^{\frac{\alpha + \beta}{2}}
(a^\dagger a - [N]_{\alpha,\beta}). \label{genhyo}
\end{subeqnarray}
Note that in the limit $q \rightarrow 1$, this algebra reduces to the
algebra for 
usual undeformed harmonic oscillator.\renewcommand{\thefootnote}{\fnsymbol{footnote}}
\footnote[4]{By considering $q^\alpha
= q^\prime$, we can actually rewrite the generalized $q$-bracket,
$[N]_{\alpha,\beta}$  as
$\dis \frac{q^\prime - q^{\prime k N }}{q - q^{\prime -1}}$, where $\dis k =
\frac{\beta}{\alpha }$. However, we find that it is sometimes more
convenient to stick to the two-parameter deformation
so that the special case of $\alpha = 0$
which corresponds effectively to $k
\rightarrow \infty$ can be discussed in similar straightforward manner
with  the case of $\beta
=0$ which corresponds to $k \rightarrow 0$.}

The general $q$-oscillator algebra is not equivalent to the Hong Yan
algebra given by eq(\ref{hyo}).  In fact, it includes the generalized
Macfarlane-Biedenharn algebra, eq(\ref{genmacf}), and  the Hong Yan
algebra, eq(\ref{hyo}), as special cases.
To see this, we first show that this general
$q$-oscillator algebra eq(\ref{genhyo})
reduces to the Macfarlane-Biedenharn
$q$-oscillator for $\alpha = 0$ or $\beta = 0$.

We note that the generalized $q$-bracket $[N]_{\alpha,\beta}$ can be
rewritten as 
\bsq
[N]_{\alpha,\beta} & = & q^{\frac{(\alpha + \beta)(N - 1)}{2}}
\frac{q^{\frac{\alpha - \beta}{2}N} - q^{- \frac{\alpha - \beta}{2}N}
}{q^{\frac{\alpha - \beta}{2}} - q^{- \frac{\alpha - \beta}{2}} } \\
& = & q^{\frac{(\alpha + \beta)(N - 1)}{2}} [N]_{q^\prime}
\esq
where $q^\prime = q^{\frac{\alpha - \beta}{2}}$.
With this result, we can recast the commutation relation in
eq(\ref{genhyo}c) as
\bsq
a a^\dagger - a^\dagger a & = & q^{\frac{(\alpha + \beta)N}{2}}
[N + 1]_{q^\prime}  - q^{\frac{(\alpha + \beta)(N - 1)}{2}}
[N]_{q^\prime} \\
q^{-\frac{(\alpha + \beta)N}{4}} a a^\dagger q^{-\frac{(\alpha +
\beta)N}{4}}  
- q^{-\frac{(\alpha + \beta)N}{2}} a^\dagger a & = &
[N + 1]_{q^\prime} - q^{\frac{-(\alpha + \beta)}{2}} [N]_{q^\prime} \\
A A^\dagger - q^{-\frac{(\alpha + \beta)}{2}} A^\dagger A & = & 
[N + 1]_{q^\prime} - q^{\frac{-(\alpha + \beta)}{2}} [N]_{q^\prime} 
\esq
Using the identity 
\beq
\lbrack N + 1 \rbrack - q \lbrack N \rbrack = q^{-N},
\eeq
one immediately
observes that the algebra in eq(\ref{genhyo}) reduces in the special
case of $\alpha = 0$ to the
Macfarlane-Biedenharn $q$-oscillator algebra with $\dis
q^\prime = q^{\frac{-
\beta}{2}}$.  Similarly, when the other parameter $\beta = 0$, 
one observes eq(\ref{genhyo}) reduces to the Macfarlane-Biedenharn
$q$-oscillator with parameter $q^{\prime \prime}= q^{\frac{-
\alpha}{2}} $. Note also that the casimir operator
$\ca_4$ in eq(\ref{genhyo}d) can be rewritten with 
the new operators $A$ and
$A^\dagger$ as $\dis q^{\frac{\alpha + \beta}{2} N} (A^\dagger
A - [N]_{\q})$.  In the limit $\alpha, \beta \rightarrow 0$, this
latter casimir
reduces to the original casimir operator $\ca_1$ in
eq(\ref{casmac1}) as we would have expected.

It is interesting to note that for $\alpha + \beta = 0$, the
general $q$-oscillator algebra eq(\ref{genhyo}) 
reduces to 
the usual Hong
Yan $q$-oscillator. Further for $\alpha =\beta$, one observes that the
generalized $q$-bracket for the operator $N$ can be rewritten as
\bsq
\lbrack N \rbrack_{\alpha,\alpha} & = & \lim_{\beta \rightarrow \alpha}
\frac{q^{\alpha N} - q^{\beta N}}{q^\alpha - q^\beta} \nonumber \\
& = & \lim_{\beta \rightarrow \alpha} q^{(N - 1)\alpha } + 
q^{(N - 2)\alpha +  \beta} + q^{(N - 3)\alpha + 2 \beta} + \cdots +
q^{(N - 1)\beta } \nonumber \\
& = & N q^{(N - 1)\alpha},
\esq   
so that the commutation relation in eq(\ref{genhyo}c) now reads
\bsq
a a^\dagger - a^\dagger a & = & q^{N \alpha} \{ 1 + N(1 - q^{-\alpha})
\} \\
\Rightarrow q^{-\frac{N}{2} \alpha} a a^\dagger q^{-\frac{N}{2} \alpha} -
q^{-\alpha} a^\dagger q^{-N \alpha} a & = & 1 + N(1 - q^{-\alpha}) \\
\Rightarrow A A^\dagger - q^{\prime \prime} A^\dagger A & = & 1 + N(1 -
q^{\prime \prime}) \label{arik}
\esq
where $q^{\prime \prime} = q^{-\alpha}$ and the operator $A$ and
$A^\dagger$ are defined by the relations $A= q^{-\frac{N}{2} \alpha} a$
and $A^\dagger = a^\dagger q^{-\frac{N}{2} \alpha}$ respectively.
Eq(\ref{arik}) can be regarded as a generalized Arik-Coon algebra.

\begin{figure}
\centerline{
\epsfxsize=80mm
\epsfysize=100mm
\epsfbox{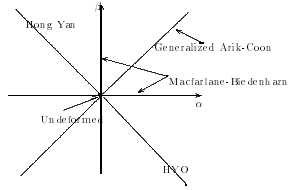}}
\vspace*{-2cm}
\caption{Schematic picture of the $\alpha-\beta$ parameter space}
\label{sche} 
\end{figure}

It is well-known that the Macfarlane-Biedenharn $q$-oscillator and the
Hong Yan $q$-oscillator cannot be mapped by invertible transformation
to each other. Since the general $q$-oscillator eq(\ref{genhyo}) 
can respectively be
reduced to these algebras for certain suitably chosen parameters, we
can conclude that the general $q$-oscillator is not
equivalent to the usual Macfarlane-Biedenharn or the Hong Yan
$q$-oscillator. Indeed, we see that in the $\alpha-\beta$ parameter
space shown in figure \ref{sche},  we have an interesting and elegant
picture in which the general $q$-oscillator algebra eq(\ref{genhyo}) seems
to interpolate between the various $q$-oscillators.
 Furthermore, if we define
$k$ as the ratio $\frac{\beta}{\alpha}$ and $\q= q^\alpha$ and
recast eq(\ref{genhyo}c) as
\beq
a a^\dagger - a^\dagger a = \frac{\q^{N+1} - \q^{k(N+1)}}{\q - \q^k}
- \frac{\q^{N} - \q^{k(N)}}{\q - \q^k}, \label{kform1}
\eeq
we can show that the general oscillator algebras
corresponding to two arbitrary $k$-values, say $k_1$ and $k_2$, are
not equivalent to each other.  To see this, one defines 
operators $\dis A = (\frac{\q - \q^{k_1}}{\q - \q^{k_2}})^{1/2}
a$ and $\dis A^\dagger = (\frac{\q - \q^{k_1}}{\q - \q^{k_2}})^{1/2}
a^\dagger$ and rewrites
eq(\ref{kform1}) for $k= k_1$ as
\beq
A A^\dagger - A^\dagger A = 
\frac{\q^{N+1} - \q^{k_2 (N+1)}}{\q - \q^{k_2}} - \frac{\q^{N} -
\q^{k_2 N}}{\q - \q^{k_2}} + {\cal F}(N, k_1, k_2)
\eeq
where $\dis {\cal F}(N,k_1, k_2) = \frac{(\q^{k_2(N+1)} - \q^{k_1(N+1)})
- (\q^{k_2 N} - \q^{k_1 N})}{\q - \q^{k_2}}$.  Clearly ${\cal F}(N, k_1,
k_2) \neq 0$ and consequently
the algebras corresponding to eq(\ref{kform1}) with
$k = k_1$ and $k = k_2$ respectively 
cannot be equivalent to each other. In
particular, one also sees that the general oscillator corresponding to
arbitrary $k$ value is not
equivalent to the Macfarlane-Biedenharn algebra which corresponds to $k = 0$ or $k
\rightarrow \infty$.

\begin{figure}
\centerline{
\epsfxsize=80mm
\epsfysize=100mm
\epsfbox{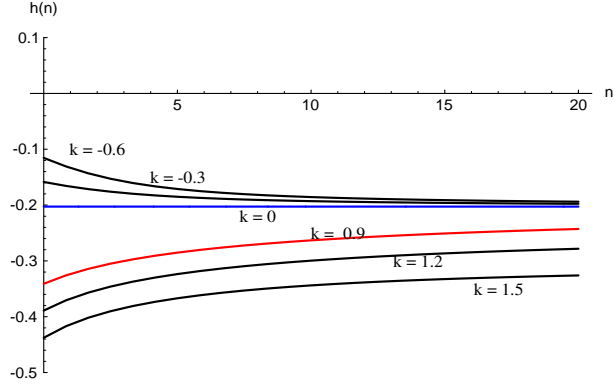}
}
\vspace*{-2cm}
\caption{Plot of $h(n)$ against $n$ for $q=1.5 $, $-0.6 \leq k \leq 1.5
$ and $\nu =0.5 $. The function $h(n)$ is independent of $n$ for $k =
0$ only.}
\label{fig2}
\end{figure}

We have already showed that the generalized Macfarlane-Biedenharn
algebra in eq (\ref{genmacf}) can be related to the usual
Macfarlane-Biedenharn oscillator, eq(\ref{macf}), by
invertible transformations of the form
$\dis A= q^{-\frac{(\alpha + \beta) N}{4}} a,
A^\dagger = a^\dagger q^{-\frac{(\alpha + \beta) N}{4}} $.  It is
interesting to explore whether there exists a  similar transformation which
can convert the general oscillator eq(\ref{genhyo})
to the Macfarlane-Biedenharn
$q$-oscillator.  To do this, one postulates existence of
new operators
\bsq B & = & f(N) a, \\
B^\dagger & = & a^\dagger f(N), \label{invert}
\esq
where $f(N)$ is an arbitrary function
of $N$ to be determined and where the
operators $\{a, a^\dagger, N \}$ and $\{B, B^\dagger, N \}$
satisfy the commutation relations for the general oscillator and the
Macfarlane-Biedenharn oscillator respectively.  Substituting
into eq(\ref{kform1}), one easily shows that
\beq
f(N)^2 B B^\dagger - f(N-1)^2 B^\dagger B = \frac{\q^{N+1} - \q^{k(N+1)}}{\q - \q^k}
- \frac{\q^{N} - \q^{k(N)}}{\q - \q^k}. \label{replace}
\eeq
Since the operators $B, B^\dagger$ are demanded to  
obey the commutation relations in
eq(\ref{macf}), one obtains the following functional
equations:
\bsq
\frac{f(N-1)^2}{f(N)^2} & = &  Q, \\
f(N)^2 (\frac{\q^{N+1} - \q^{k(N+1)}}{\q - \q^k}
- \frac{\q^{N} - \q^{k(N)}}{\q - \q^k} ) & = &  Q^{-N}, \label{func}
\esq
where $Q$ is some function of $\q$ and $k$ and should be independent of
the operator $N$. 
Solving eq(\ref{func}a) inductively gives $\dis f(N) = Q^{-\frac{N}{2}}$
which leads to the relation
\beq
\frac{\q^{N+1} - \q^{k(N+1)}}{\q - \q^k}
- \frac{\q^{N} - \q^{k(N)}}{\q - \q^k}   =  Q^{-2N}. \label{solvQ}
\eeq
One needs to solve eq(\ref{solvQ})
for $Q$.  However it can be shown, at least at the representation level,
that this is not possible except for
$k = 0$, which corresponds to the Macfarlane-Biedenharn case.  
More specifically, one can always specify
a representation $\{ |n> \}$ in which $N |n> = (n + \nu) |n>$ as in ref
\cite{rideau} and investigate the function $\ln Q \equiv h(n)$,
$$
h(n) = \frac{-1}{2(n + \nu)}\ln \left( \frac{\q^{n + \nu +1} 
- \q^{k( n + \nu +1)}}{\q - \q^k}
- \frac{\q^{n + \nu} - \q^{k(n + \nu)}}{\q - \q^k} \right), 
$$ to see whether it is $n$-independent. 
Here, $\nu$ is a free parameter.
Sweeping over a wide range values for $k$ and $\nu$,
we find graphically that
$h(n)$ does depend on $n$.  A plot of this function for typical values of
$q, k$ and $\nu$ is displayed in figure \ref{fig2} and shows
clearly that $h(n)$ depends on $n$ except for $k = 0$.
Moreover, using graphical method we find 
the partial derivative $\frac{\partial h}{\partial n}$ 
does not vanish unless $k=0$.
Thus, one concludes that it is not possible
to relate the general oscillator eq(\ref{genhyo}) 
for any arbitrary $k$ (except $k = 0$) to the
Macfarlane-Biedenharn $q$-oscillator eq(\ref{macf}) 
through the invertible
transformations in eq(\ref{invert}).

Finally, we note that since the Hong Yan $q$-oscillator can be endowed
with a Hopf structure, it would be natural to question the existence of
a similar Hopf structure for the general $q$-oscillator.  
Our preliminary work seems to indicate that a Hopf structure for the
general oscillator exists only for the special case of $\alpha = -
\beta$. There are several reasons for this suspicion.  Firstly,
one can postulate the existence of the following coproducts \cite{oh,tso}
\bsq
\Delta(a^\dagger) &=& c_1 a^\dagger \otimes \q^{\alpha_1 N} + c_2
\q^{\alpha_2 N} \otimes a^\dagger \\
\Delta(a) & = & c_3 a \otimes \q^{\alpha_3 N} + c_4 \q^{\alpha_4 N}
\otimes a \\
\Delta(N) & = & c_5 N \otimes 1 + c_6 1 \otimes N + \gamma 1 \otimes 1 \label{asco}
\esq
where $c_i (i = 1,2, \cdots 6)$, $\alpha_i (i=1,2, \cdots 4)$ and
$\gamma$ are constants to be determined. One then demands that the
coproducts satisfy the
co-associativity axiom, namely 
\beq
(id \otimes \Delta)\Delta(h) = (\Delta \otimes id) \Delta(h),
\eeq
where $h$ is any of the generators $\{ a, a^\dagger, N \}$. Furthermore,
one must also ensure that the following consistency condition
\beq
\Delta(a) \Delta(a^\dagger) - \Delta(a^\dagger)\Delta(a) =
\Delta(\frac{\q^{N+1} - \q^{k(N + 1)}}{\q - \q^k}) -
\Delta(\frac{\q^{N} - \q^{kN}}{\q - \q^k}) 
\eeq
is satisfied.
With these constraints, a
straightforward calculation shows that $k$ must necessarily be $-1$ 
which indicates that the coproducts in eq(\ref{asco}) can only support 
the Hong Yan
oscillator.
Secondly, Hong Yan algebra can be truncated and
identified with $su_{\sqrt{q}} (2)$ algebra for some values of $q$.  
With this identification, one then easily borrows
the Hopf structure for $su_{\sqrt{q}} (2)$ algebra for the Hong Yan
oscillator algebra.  A similar
truncation for the general oscillator does not seem possible.  
Lastly,  a naive
generalization of $su_{q} (2)$ algebra in the form
\bsq
\lbrack J_0, J_\pm \rbrack & = & \pm J_\pm \\
\lbrack J_+, J_- \rbrack & = & \frac{q^{2 \alpha J_0} - q^{2 \beta
J_0}}{q^\alpha - q^\beta} = [2 J_0]_{\alpha,\beta}
\esq
does not lead to a generalized coproduct which can maintain the algebra
homomorphism $\Delta [J_+,
J_-] = [\Delta(J_+), \Delta(J_-)]$, 
unless we set $\alpha = - \beta$. However,
it is still premature to say conclusively that the Hong
Yan limit is the only one in which one can 
associate the general oscillator with a Hopf
structure. 
Further development in this
direction would probably help us gain invaluable insights into the
nature of the $q$-deformation of the harmonic oscillator.


\end{document}